\documentclass[11pt]{article}

\usepackage{amsmath,amsfonts,amssymb} 

\usepackage{amsthm} 

\usepackage{graphicx,caption,subcaption} 
\usepackage{fullpage}

\allowdisplaybreaks 

\begingroup
    \makeatletter
    \@for\theoremstyle:=definition,remark,plain\do{%
        \expandafter\g@addto@macro\csname th@\theoremstyle\endcsname{%
            \addtolength\thm@preskip\parskip
            }%
        }
\endgroup

\newtheorem{thm}{Theorem}
\newtheorem{cor}[thm]{Corollary}

\newtheorem{lemma}[thm]{Lemma}

\begin{document}

\title{\vspace{-0.5in} On a property of $2$-connected graphs and Dirac's Theorem}

\author{
{{Alexandr Kostochka}}\thanks{
\footnotesize {University of Illinois at Urbana--Champaign, Urbana, IL 61801. E-mail: \texttt {kostochk@math.uiuc.edu}.
 Research 
is supported in part by  NSF  Grant DMS-2153507.}}
\and
{{Ruth Luo}}\thanks{
\footnotesize {University of South Carolina, Columbia, SC 29208, USA. E-mail: \texttt {ruthluo@sc.edu}.
}}
\and{{Grace McCourt}}\thanks{University of Illinois at Urbana--Champaign, Urbana, IL 61801, USA. E-mail: {\tt mccourt4@illinois.edu}. Research 
is supported in part by NSF RTG grant DMS-1937241.}}

\date{ \today}
\maketitle

\vspace{-0.3in}

\begin{abstract}
We refine a property of $2$-connected graphs described in the classical paper of Dirac from 1952 and use the refined property to somewhat shorten Dirac's proof of the fact that each $2$-connected $n$-vertex graph with minimum degree at least $k$ has a cycle of length at least $\min\{n,2k\}$.

\end{abstract}

\section{Introduction}
One of the basic facts on $2$-connected graphs is their characterization by Whitney~\cite{Wh} from 1932:

\begin{thm}[Whitney~\cite{Wh}]\label{wh} 
A  graph $G$ with at least $3$ vertices is $2$-connected if and only if for any distinct $u,v\in V(G)$ there are internally disjoint $u,v$-paths.
\end{thm}

Given  paths $P$ and $P'$ with the common origin  in a graph, we say $P'$ is {\em aligned with} $P$ if for all $u,v \in V(P) \cap V(P')$ if $u$ appears before $v$ in $P$, then $u$ also appears before $v$ in $P'$. 

In his thesis~\cite{diracPh} and classical paper~\cite{dirac}, Dirac refined (the main part of)  Theorem~\ref{wh} as follows.

\begin{lemma}[Dirac, Lemma 2 in~\cite{dirac}]\label{dirl2} 
If $x$ and $y$ are two distinct vertices of a graph without cut vertices, and if $W$ is any given path connecting $x$ and $y$, then the graph contains two paths connecting $x$ and $y$ and having the following properties: (i) they are internally disjoint;
(ii) each of them is aligned with $W$.
\end{lemma}

He also says that this lemma can be further refined as follows.

\begin{cor}[Dirac, Corollary on p.73 in~\cite{dirac}]\label{dirl21}  
If $x$ is adjacent to a vertex  $z$ of $W$, then the graph contains two
paths connecting $x$ and $y$ such that they have the properties (i) and (ii), and
one of them goes through $z$.
\end{cor}

Dirac used this corollary to prove the following famous theorem:

\begin{thm}[Dirac, Theorem 4 in~\cite{dirac}]\label{dirt} 
Let $n>k\geq 2$. Each $2$-connected $n$-vertex graph with minimum degree at least $k$ has a cycle of length at least $\min\{n,2k\}$.
\end{thm}
P\' osa~\cite{Po} used Lemma~\ref{dirl2} to derive an extension of Theorem~\ref{dirt}.
In this note, we refine  Corollary~\ref{dirl21} (with an almost the same proof) as follows.

\begin{lemma}\label{ndirl2} 
Let $P$ be an $x,y$-path in a 2-connected graph $G$, and let $z \in V(G)$. Then there exists an $x, z$-path $P_1$ and an $x,y$-path $P_2$ such that

(a) $P_1$ and $P_2$ are internally disjoint, and
$\quad$ (b) each of $P_1$ and $P_2$ is aligned with $P$.
\end{lemma}

We then use this lemma to give a somewhat shorter and logically simpler proof of Theorem~\ref{dirt}.

\medskip
{\bf Remark 1.} The authors~\cite{KLM} 
used Lemma 5 to prove an analog of Theorem 4 for Berge cycles in $r$-uniform 2-connected hypergraphs.

{\bf Remark 2.} Douglas West pointed out how to easily  derive Lemma~\ref{ndirl2} from Lemma~\ref{dirl2}. We still present the full proof in order to show a shorter proof from scratch for Theorem~\ref{dirt}.

{\bf Remark 3.} The straightforward generalization of Lemma~\ref{ndirl2} or Lemma~\ref{dirl2} to $k$-connected graphs is not true. In fact, for {\em each} positive integer $k$, there exists a $k$-connected graph $G$ and an $x,y$-path $P$  such that $G$ has no $3$  $x,y$-paths aligned with $P$. 

\begin{figure}[h]
    \centering
    \begin{minipage}{0.43\textwidth}
        \centering
        \includegraphics[width=0.7\textwidth]{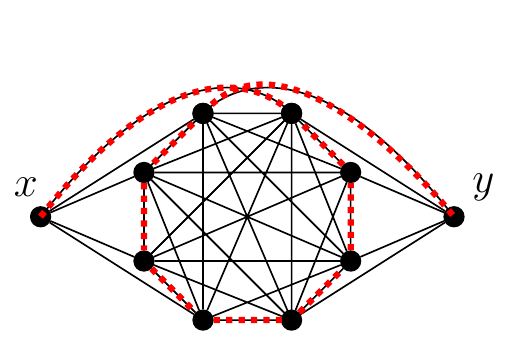} 
        \caption{A 5-connected construction}
        \label{aligned3}
        
    \end{minipage}
    \begin{minipage}{0.5\textwidth}
        \centering
        \includegraphics[width=1.0\textwidth]{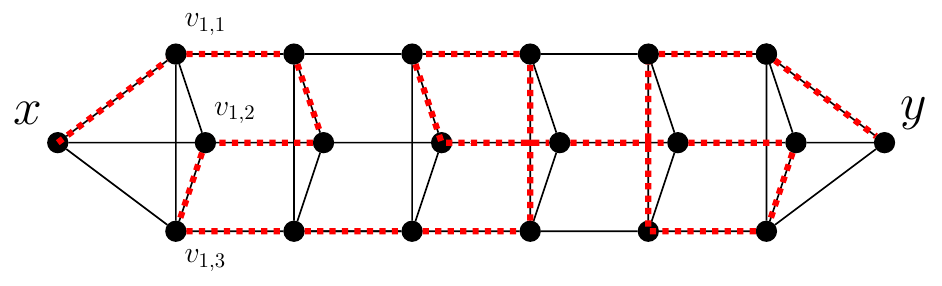} 
        \caption{Construction $H_3$}
        \label{alignedpath1}
    \end{minipage}
\end{figure}

Set $A = \{a_1, \ldots, a_{k-1}\}$ and $B =\{b_1, \ldots, b_{k-1}\}$, and let $G_k$ be the graph with vertex set $A \cup B \cup \{x,y\}$ such that $G_k[A \cup B]$ induces a clique, $N(x) = A \cup \{b_1\}$ and $N(y) = B \cup \{a_1\}$. Define the path $P=x, b_1, b_2, \ldots, b_{k-1}, a_{k-1}, a_{k-2}, \ldots, a_1, y$. Observe that $G_k$ is $k$-connected and any $x,y$-path aligned in $P$ in $G$ must use  edge $xb_1$ or edge $a_1y$. Hence $G_k$ has at most $2$  internally disjoint paths aligned with $P$. Fig.~\ref{aligned3} displays the construction for $k=5$.


{\bf Remark 4.} Let $x$ and $y$ be distinct vertices in a graph $G$. If $G$ is $2$-connected, then
 Lemma 2 provides two internally disjoint $x,y$-paths aligned with any  fixed $x,y$-path $P$ in $G$. If  $G$ is $k$-connected for some $k\geq 3$, then Menger's Theorem guarantees $k$  internally disjoint $x,y$-paths. A natural question is: Given $k\geq 3$, a $k$-connected graph $G$ and
  an $x,y$-path $P$ in $G$, how many of the $k$ paths in Menger's Theorem always can be chosen  to be aligned with $P$?
  
 Somewhat surprisingly, the answer is ``zero". We construct an example $H_3$ for $k=3$ is as follows (see Fig.~\ref{alignedpath1}). 
 For $1\leq i\leq 6$, let $W_i=\{v_{i,1},v_{i,2},v_{i,3}\}$. Let $V(H_3)=\{x,y\}\cup \bigcup_{i=1}^{6}W_i$. For $1\leq j\leq 3$, let $P_j=x,v_{1,j},v_{2,j},\ldots,v_{6,j},y$. The edge set of $H_3$ contains the edges of these three paths plus $H_3[W_i]=K_3$ for all $1\leq i\leq 6$. By construction, $H_3$ is $3$-connected, and
  the only triple of internally disjoint $x,y$-paths in $H_3$ is $\{P_1,P_2,P_3\}$. But none of these paths is aligned with the $x,y$-path
 $$P_0=x,v_{1,1},v_{2,1},v_{2,2},v_{1,2},v_{1,3},v_{2,3},v_{3,3},v_{4,3},v_{4,1},v_{3,1},v_{3,2},v_{4,2},
v_{5,2},v_{6,2},v_{6,3},v_{5,3},v_{5,1},v_{6,1},y:
 $$
the edge $v_{4,1}v_{3,1}$ in $P_0$ is opposite to an edge in $P_1$,
the edge $v_{2,2}v_{1,2}$  is opposite to an edge in $P_2$ and the edge $v_{6,3}v_{5,3}$  is opposite to an edge in $P_3$.

The examples for $k\geq 4$ are very similar.

{\bf Remark 5.} In the statement of Lemma~\ref{ndirl2}, we cannot replace an $x,y$-path $P_2$ with an $x,z'$-path for an arbitrary $z'$: If $n\geq 5$, $G$ is an $n$-cycle $v_1,v_2,\ldots,v_n,v_1$ and $P=v_1,v_2,\ldots,v_n$, then $G$ has no two $v_1,v_3$-paths both aligned with $P$.

{\bf Remark 6.} After the proof of Theorem~\ref{dirt}, Dirac~\cite{dirac} thanks a referee for simplifying his original proof. Bjarne Toft\footnote{Private communication.} suggests that this referee possibly was Harold A. Stone.

\section{Proofs}

We view paths as having one of the two possible orientations.
For a path $P$ and $u,v,w\in V(P)$, let $P[u,v]$ denote the subpath of $P$ from $u$ to $v$, and let $P^+(w)$ denote the part of $P$ starting from $w$.

 \subsection{Proof of   Lemma~\ref{ndirl2} }

We induct on $|V(P)|$. The base case $P = x,y$ follows from the fact that the connected graph $G-xy$ has an $x,y$-path $P_1$, so we can take $P_2=P$. 

Induction step: Let the lemma hold for all paths with fewer than $s$ vertices and let 
$P = v_1, v_2, \ldots, v_s$ with $x=v_1$ and $y=v_s$. By the induction hypothesis, the lemma holds for $P'=v_2,v_3,\ldots,v_s$ with $x'=v_2$ and $y'=v_s$. Let 
  $P'_1$ and $P'_2$ be the corresponding $v_2,z$-path and $v_2,v_s$-path, respectively.

{\bf Case 1}.  $v_1\in P'_i$ for some $i \in \{1,2\}$. Then let $P_i=(P'_i)^+(v_1)$, and let $P_{3-i}$ be obtained  $P'_{3-i}$ by adding  edge $v_{1}v_2$ at the start.

{\bf Case 2}. $v_1 \notin P'_i$ for $i =1,2$. 
Since $G$ is $2$-connected, $G-v_{2}$ has a path from $v_1$ to $P\cup P'_1\cup P'_2$.
Let $Q$ be a shortest such path and $u$ be the end of $Q$ distinct from $v_1$.

If $u\in P'_i$, then replace $P'_i$ with $Q,(P'_i)^+(u)$, and then extend $P'_{3-i}$ by adding  $v_{1}v_2$.
Suppose now that
 $ u\in P-P'_1-P'_2$, say $u=v_j$. Choose the minimum $j'\geq j$ such that $v_{j'} \in P'_1\cup P'_2$. It exists, since $v_s\in P_2$, say 
$v_{j'} \in P'_i$.
 In this case,   let $P_i=Q,P[v_{j},v_{j'}],(P'_i)^+(v_{j'})$, and  $P_{3-i}=v_1,v_2,P'_{3-i}$. This proves the lemma.

 \subsection{Proof of Theorem~\ref{dirt}}
 
Suppose graph $G$ is a counter-example to the theorem and its minimum degree, $\delta(G)$, is $k$. Since $G$ is $2$-connected, $k\geq 2$.

A {\em lollipop} in $G$ is a pair $(C,P)$ where $C$ is a cycle and $P$ is a path such that $V(C)\cap V(P)$ is one vertex that is an end of $P$ (see Fig.~\ref{lollipop} below). A lollipop $(C,P)$ is {\em better than} a lollipop $(C',P')$ if $|V(C)|>|V(C')|$ or
 $|V(C)|=|V(C')|$ and  $|V(P)|>|V(P')|$.

\begin{figure}[h]
        \centering
        \includegraphics[width=0.35\textwidth]{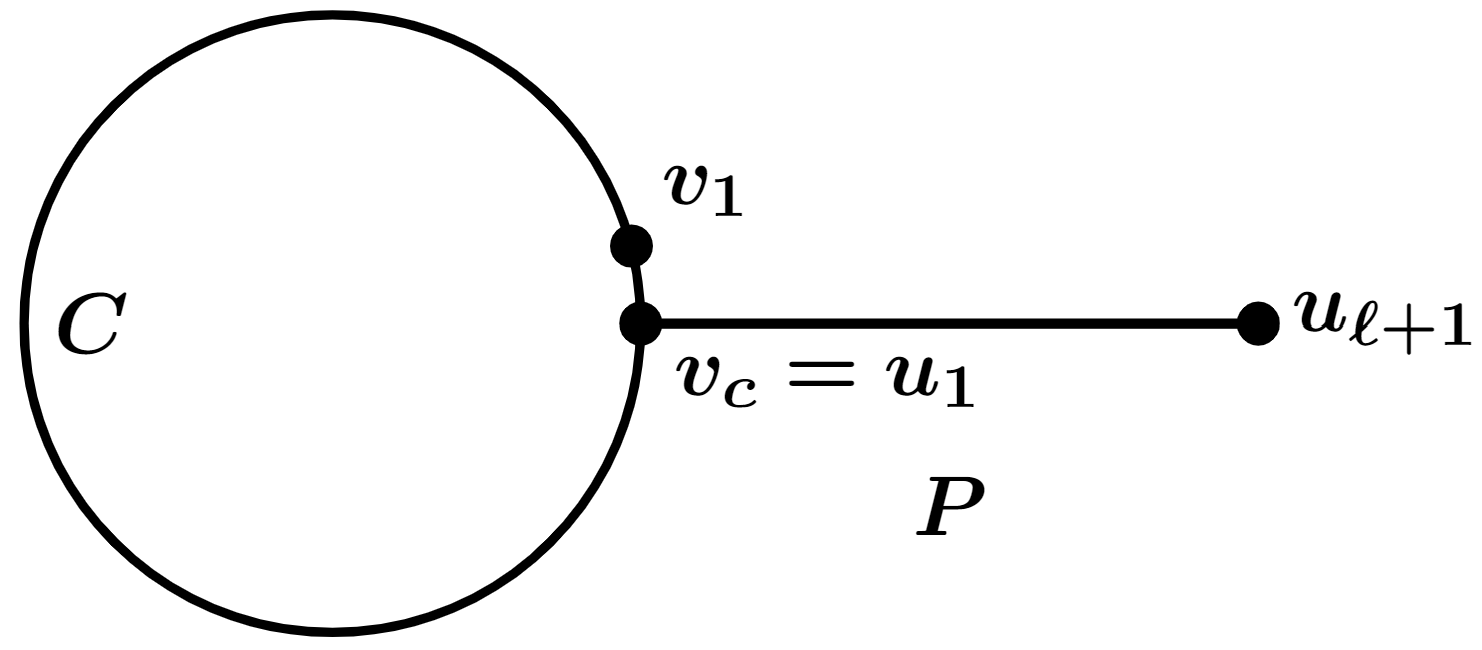} 
        \caption{A lollipop $(C,P)$}
        \label{lollipop}
\end{figure}

 Let $(C,P)$ be a best lollipop in $G$. For definiteness, let $C=v_1,v_2,\ldots,v_c,v_1$ and
 $P=u_1,u_2,\ldots,u_{\ell+1}$, where $u_1=v_c$. Since $G$ is a counterexample, $c<n$. So, since $G$ is
 $2$-connected, $\ell\geq 1$.

 By the maximality of $\ell$, $N(u_{\ell+1})\subseteq V(P)\cup V(C)$.
 
 \medskip
 {\bf Case 1:} There is $v_i\in N(u_{\ell+1})$. If $1\leq i\leq \ell$, then the cycle $v_i,v_{i+1},\ldots,v_c,u_2,u_3,\ldots,u_{\ell+1},v_i$ is longer than $C$, a contradiction. Thus $i\geq \ell+1$. Symmetrically,
 $i\leq c-\ell-1$. On the other hand, if $N(u_{\ell+1})\cap V(C)$ contains two consecutive vertices $v_i$ and $v_{i+1}$, then replacing edge $v_iv_{i+1}$ in $C$ with the path $v_i,u_{\ell+1},v_{i+1}$, we again get a cycle longer than $C$. 
 Since $|N(u_{\ell+1})\cap V(P)|\leq \ell$, we get
 $$k\leq d(u_{\ell+1})=|N(u_{\ell+1})\cap V(P)|+|N(u_{\ell+1})\cap V(C)-v_c|\leq \ell+\left\lceil \frac{c-1-2\ell}{2}\right\rceil
 =\left\lceil \frac{c-1}{2}\right\rceil<k,$$
 a contradiction.

 \medskip
 {\bf Case 2:} $N(u_{\ell+1})\subseteq V(P)$. Let $N(u_{\ell+1})=\{u_{j_1},\ldots,u_{j_s}\}$ with $j_1<j_2<\ldots,<u_s$.
 Let $P'=v_1,v_2,\ldots,v_c,u_2,\ldots,u_{\ell+1}$, $x=v_1$, $y=u_{\ell+1}$ and $z=u_{j_1}$. By Lemma~\ref{ndirl2} for these $P',x,y$ and $z$, there exists an $x,z$-path $P_1$ and an $x,y$-path $P_2$ that are internally disjoint and aligned with $P'$.
 
 For $h\in \{1,2\}$, let $a_h$ be the last vertex of $P_h$ in $C$ and $b_h$ be the first vertex of $P_h$ in 
 $Y=\{u_{j_1},u_{j_1+1},\ldots,u_{\ell+1}\}$. Since $P_1$ is aligned with $P'$, $b_1=u_{j_1}$.
 
 If $a_2=a_1$, then since $P_1$ and $P_2$ are internally disjoint, $a_2=a_1=v_1$, and one of $P_1$ and $P_2$, say $P_h$, does not contain $v_c (=u_1)$ and first intersects $P$ at some vertex $u_i$ with $i \geq 2$. Then deleting from $C$ edge $v_1v_c$ and adding instead paths $P_h[v_1,u_i]$ and $P[u_1,u_i]$, we get a longer cycle, a contradiction. Thus, $P_1[a_1,b_1]$ and $P_2[a_2,b_2]$ are disjoint.
 
 Let $Q$ be the longer of the two subpaths of $C$ connecting $a_1$ with $a_2$. Then $|V(Q)|\geq 1+\left\lceil \frac{c}{2}\right\rceil$. Let $b_2=u_j$. Since $P_1$ and $P_2$ are internally disjoint, $j>j_1$.
  Let $j'$ be the largest index in $\{j_1,\ldots,j_s\}$ that is less than $j$. Since $j>j_1$, $j'$ is well defined and $j'\geq j_1$.
  
Consider the closed walk 
$$C'=a_2,Q,a_1,P_1[a_1,u_{j_1}],u_{j_1},u_{j_1+1},\ldots,u_{j'},u_{\ell+1},u_{\ell},\ldots,u_{j}(=b_2),P_2[b_2,a_2],a_2.$$
Since $P_1[a_1,u_{j_1}]$ and $P_2[b_2,a_2]$ are disjoint and both are internally disjoint from $V(C)\cup Y$, $C'$ is a cycle.
It has at least $|V(Q)|\geq 1+\left\lceil \frac{c}{2}\right\rceil=1+c-\left\lfloor\frac{c}{2}\right\rfloor$ vertices in $C$ and at least $1+d(u_{\ell+1})\geq 1+k$ vertices in $P$. Since $|V(C)\cap V(P)|=1$, we have 
$$|V(C')|\geq (1+c-\left\lfloor\frac{c}{2}\right\rfloor)+(1+k)-1>c.
$$
This contradiction with the maximality of $C$ proves the theorem.

\bigskip
{\bf Acknowledgment.} We thank Bjarne Toft and Douglas West for their valuable comments.

\end{document}